\documentclass[10pt]{amsart}
\usepackage{graphicx,amsmath,amsthm,amssymb,verbatim,enumerate,ifthen}
\usepackage[utf8]{inputenc}
\def\N{\mathbb{N}}
\def\R{\mathbb{R}}

\def\conv{\mathop{\mbox{\rm conv}}}
\def\comment#1{}

\newtheorem{theorem}{Theorem}
\def\Thm#1#2{\begin{theorem}\label{T#1}#2\end{theorem}}
\def\thm#1{Theorem~\ref{T#1}}

\newtheorem{corollary}[theorem]{Corollary}
\def\Cor#1#2{\begin{corollary}\label{C#1}#2\end{corollary}}
\def\cor#1{Corollary~\ref{C#1}}

\newtheorem{lemma}[theorem]{Lemma}
\def\Lem#1#2{\begin{lemma}\label{L#1}#2\end{lemma}}
\def\lem#1{Lemma~\ref{L#1}}

\theoremstyle{definition}
\newtheorem{definition}[theorem]{Definition}
\def\Dfn#1#2{\begin{definition}\label{Defi#1}\rm #2\end{definition}}

\def\eq#1{{\rm(\ref{E#1})}}
\def\Eq#1#2{\ifthenelse{\equal{#1}{*}}
  {\begin{equation*}\begin{aligned}#2\end{aligned}\end{equation*}}
  {\begin{equation}\begin{aligned}\label{E#1}#2\end{aligned}\end{equation}}}

\begin{document}
\date{\today}
\title{A maximum theorem for generalized convex functions}

\author[Zs. P\'ales]{Zsolt P\'ales}
\address{Institute of Mathematics, University of Debrecen, H-4002 Debrecen, Pf.\ 400, Hungary}
\email{pales@science.unideb.hu}


\dedicatory{Dedicated to the memory of Professors Gábor Kassay and Csaba Varga}

\thanks{The research of the author was supported by the K-134191 NKFIH Grant and the 2019-2.1.11-TÉT-2019-00049 project. }

\begin{abstract}
Motivated by the Maximum Theorem for convex functions (in the setting of linear spaces) and for subadditive functions (in the setting of Abelian semigroups), we establish a Maximum Theorem for the class of generalized convex functions, i.e., for functions $f:X\to X$ that satisfy the inequality $f(x\circ y)\leq pf(x)+qf(y)$, where $\circ$ is a binary operation on $X$ and $p,q$ are positive constants. As an application, we also obtain an extension of the Karush--Kuhn--Tucker theorem for this class of functions.
\end{abstract}

\maketitle

\section{Introduction}


In what follows, a linear space $X$ always means a vector space over the field of real numbers. If $X$ is a topological linear space, then its (topological) dual space is denoted by $X^*$. The Maximum Theorem for convex functions, which is due to Dubovitskii and Milyutin (cf.\ \cite{MagTik03}), can be stated as follows. 

\Thm{MTD}{Let $X$ be a linear space, let $D\subseteq X$ be a convex set and let $f_1,\dots,f_n:D\to\R$ be convex functions such that
\Eq{*}{
  0\leq\max(f_1(x),\dots,f_n(x)) \qquad (x\in D).
}
Then there exist $\lambda_1,\dots,\lambda_n\geq0$ with
$\lambda_1+\cdots+\lambda_n=1$ such that
\Eq{*}{
  0\leq \lambda_1f_1(x)+\cdots+\lambda_nf_n(x) \qquad (x\in D).
}}

A standard application of the Maximum Theorem is to prove the subdifferential formula for the pointwise maximum of convex functions, which was established by Dubovitskii and Milyutin (see \cite{MagTik03}). For the standard terminologies and notations, we refer to the list of monographs in the list of references, where the reader can find many more details and applications.

\Thm{SD}{Let $X$ be a topological vector space, $D\subseteq X$ be an open convex set, $p\in D$ and $f_1,\dots,f_n:D\to\R$ be continuous convex functions with $f_1(p)=\cdots=f_n(p)$ and define $f:=\max(f_1,\dots,f_n)$. Then
\Eq{*}{
  \partial f(p)
  =\conv\big(\partial f_1(p)\cup\dots\cup\partial f_n(p)\big).
}}

\begin{proof}
Using that $f(p)=f_1(p)=\cdots=f_n(p)$, for all $h\in X$, we obtain
\Eq{*}{
   f'(p,h):&=\lim_{t\to0^+}\frac{f(p+th)-f(p)}{t}\\
   &=\lim_{t\to0^+}\frac{\max(f_1(p+th),\dots,f_n(p+th))-f(p)}{t}\\
   &=\lim_{t\to0^+}\max\bigg(\frac{f_1(p+th)-f(p)}{t},\dots,\frac{f_n(p+th)-f(p)}{t}\bigg)\\
   &=\lim_{t\to0^+}\max\bigg(\frac{f_1(p+th)-f_1(p)}{t},\dots,\frac{f_n(p+th)-f_n(p)}{t}\bigg)\\
   &=\max\bigg(\lim_{t\to0^+}\frac{f_1(p+th)-f_1(p)}{t},\dots,\lim_{t\to0^+}\frac{f_n(p+th)-f_n(p)}{t}\bigg)\\
   &=\max(f'_1(p,h),\dots,f'_n(p,h)).
}
First assume that a continuous linear functional $\varphi \in X^*$ belongs to $\partial f(p)$. Then, in view of the above formula for directional derivatives, we get 
\Eq{*}{
  \varphi(h)\leq f'(p,h)=\max(f'_1(p,h),\dots,f'_n(p,h))\qquad(h\in X).
}
This relation implies that
\Eq{*}{
  0\leq \max(f'_1(p,h)-\varphi(h),\dots,f'_n(p,h)-\varphi(h))\qquad(h\in X).
}
This inequality states that the maximum of the convex functions $h\mapsto f'_i(p,h)-\varphi(h)$ is nonnegative.
Thus, by the Maximum Theorem, there exist $\lambda_1,\dots,\lambda_n\geq0$ with $\lambda_1+\cdots+\lambda_n=1$ such that 
\Eq{*}{
  0\leq \lambda_1(f'_1(p,h)-\varphi(h))+\cdots+\lambda_n(f'_n(p,h)-\varphi(h))\qquad(h\in X),
}
equivalently,
\Eq{*}{
  \varphi(h)
  \leq \lambda_1f'_1(p,h)+\cdots\lambda_nf'_n(p,h)
  =(\lambda_1f_1+\cdots\lambda_nf_n)'(p,h)
  \qquad(h\in X).
}
Using the so-called Sum Rule, we get
\Eq{*}{
  \varphi\in\partial(\lambda_1f_1+\cdots\lambda_nf_n)(p)
  &=\lambda_1\partial f_1(p)+\cdots+\lambda_n\partial f_n(p)\\
  &\subseteq\conv\big(\partial f_1(p)\cup\dots\cup\partial f_n(p)\big).
}
The proof of the reversed inclusion is simpler, thus it is left to the reader.
\end{proof}


Another motivation for this paper comes from the theory of subadditive functions defined on Abelian semigroups. The following result was stated in the monograph \cite{FucLus81} of Fuchssteiner and Lusky.

\Thm{MTSD}{Let $(X,+)$ be an Abelian semigroup and let $f_1,\dots,f_n:X\to\R$ be subadditive functions such that
\Eq{*}{
  0\leq\max(f_1(x),\dots,f_n(x)) \qquad (x\in X).
}
Then there exist $\lambda_1,\dots,\lambda_n\geq0$ with
$\lambda_1+\cdots+\lambda_n=1$ such that
\Eq{*}{
  0\leq \lambda_1f_1(x)+\cdots+\lambda_nf_n(x) \qquad (x\in X).
}}

This result has beautiful applications in the book \cite{FucLus81}, for instance, the Phragmen--Lindel\"of Principle and the Hadamard Three Circle Theorem (both results belong to the theory of complex functions) can elegantly be verified in terms of them. 

\section{The General Maximum Problem}

The two Maximum Theorems described in the Introduction motivate the following definition.

\Dfn{1}{Let $X$ be a nonempty set. A family $\mathcal{F}\subseteq\{f:X\to\R\}$ is said to have the \emph{discrete maximum property} if 
\Eq{*}{
  f_1,\dots,f_n\in\mathcal{F},\qquad 0\leq\max(f_1(x),\dots,f_n(x))\qquad(x\in X)
}
implies that there exist $(\lambda_1,\dots,\lambda_n)\in S_n$ such that
\Eq{*}{
  0\leq \lambda_1f_1(x)+\cdots+\lambda_nf_n(x)\qquad(x\in X).
}
Here, for convenience, $S_n$ denotes the $(n-1)$-dimensional simplex
\Eq{*}{
  \{(\lambda_1,\dots,\lambda_n)\in\R^n\mid \lambda_1,\dots,\lambda_n\geq0,\,\lambda_1+\dots+\lambda_n=1\}.
}}

If $X$ has at least two elements, then the set of all functions $\mathcal{F}:=\{f:X\to\R\}$ does not have the discrete maximum property. Indeed, Let $\{A_1,A_2\}$ be a partition of $X$ and $f_i(x):=0$ if $x\in A_i$, $f_i(x):=-1$ if $x\not\in A_i$. Then $\max(f_1,f_2)=0$, but $\lambda f_1+(1-\lambda)f_2<0$ for all $\lambda\in[0,1]$. This example shows that, in order to possess the discrete maximum property, the family $\mathcal{F}\subseteq\{f:X\to\R\}$ must satisfy some additional nontrivial conditions. 

In the next result we characterize the situation when a finite family of given functions possess a nonnegative convex combination. 

\Thm{nF}{Let $X$ be nonempty and $f_1,\dots, f_n:X\to\R$. Then there exists $(\lambda_1,\dots,\lambda_n)\in S^n$ such that 
\Eq{lfx}{
  0\leq\lambda_1f_1(x)+\dots+\lambda_nf_n(x)\qquad(x\in X)
}
if and only if
\Eq{xt}{
  0\leq \max_{i\in\{1,\dots,n\}}
  \big(t_1f_i(x_1)+\dots+t_nf_i(x_n)\big)
  \quad(x_1,\dots,x_n\in X,\,(t_1,\dots,t_n)\in S_n).
}}

\begin{proof}
Assume first that \eq{lfx} holds for some $\lambda\in S_n$. To verify the necessity of \eq{xt}, let $x_1,\dots,x_n\in X$ and $(t_1,\dots,t_n)\in S_n$ be arbitrary. Then, using \eq{lfx} for $x\in\{x_1,\dots,x_n\}$, we get
\Eq{*}{
  0&\leq \sum_{j=1}^n t_j\big(\lambda_1f_1(x_j)+\dots+\lambda_nf_n(x_j)\big)\\
  &=\sum_{i=1}^n \lambda_i\big(t_1f_i(x_1)+\dots+t_nf_i(x_n)\big)\\
  &\leq \max_{i\in\{1,\dots,n\}}
  \big(t_1f_i(x_1)+\dots+t_nf_i(x_n)\big).
}
This shows the necessity of condition \eq{xt}.

Now assume that \eq{xt} holds and, for $x\in X$, define
the set $\Lambda_x\subseteq S_n$ by
\Eq{Lx}{
  \Lambda_x:=\big\{(\lambda_1,\dots,\lambda_n)\in S_n\mid 
  0\leq\lambda_1f_1(x)+\dots+\lambda_nf_n(x)\big\}.
}
The inequality \eq{lfx} is now equivalent to the condition
\Eq{IS}{
  \bigcap_{x\in X}\Lambda_x\neq\emptyset,
}
because every element $\lambda$ of the above intersection will satisfy \eq{lfx}. It easily follows from the definition that $\Lambda_x$ is a compact convex subset of the $(n-1)$-dimensional affine space 
\Eq{*}{
  \{(\lambda_1,\dots,\lambda_n)\in\R^n\mid \lambda_1+\dots+\lambda_n=1\}.
}
Therefore, according to Helly's Theorem, the condition \eq{IS} is satisfied if and only every $n$-member subfamily of $\{\Lambda_x\mid x\in X\}$ has a nonempty intersection. To verify this, let $x_1,\dots.x_n\in X$ be fixed arbitrarily. According to inequality \eq{xt}, the pointwise maximum of the convex functions
\Eq{*}{
  S_n\ni(t_1,\dots,t_n)
  \mapsto t_1f_i(x_1)+\dots+t_nf_i(x_n)
}
is nonnegative over $S_n$. Therefore, in view of \thm{MTD}, there exists $(\lambda_1,\dots,\lambda_n)\in S_n$ such that 
\Eq{*}{
  0&\leq\sum_{i=1}^n\lambda_i\big(t_1f_i(x_1)+\dots+t_nf_i(x_n)\big)\\
  &=\sum_{j=1}^n t_j\big(\lambda_1f_1(x_j)+\dots+\lambda_nf_n(x_j)\big) \qquad ((t_1,\dots,t_n)\in S_n).
}
If $i\in\{1,\dots,n\}$, then substituting $(t_1,\dots,t_n):=(\delta_{i,j})_{j=1}^n$ into the above inequality, we get that
\Eq{*}{
  \lambda_1f_1(x_i)+\dots+\lambda_nf_n(x_i) \qquad(i\in\{1,\dots,n\}).
}
This shows that $\lambda\in\Lambda_{x_1}\cap\dots\cap\Lambda_{x_n}$, proving that this intersection is nonempty, as it was desired.
\end{proof}

In the case $n=2$, the above theorem immediately implies the following statement.

\Cor{2F}{Let $X$ be a nonempty set and $f,g:X\to\R$. Then there exists $\lambda\in[0,1]$ such that 
\Eq{fgx}{
  0\leq \lambda f(x)+(1-\lambda) g(x) \qquad(x\in X)
}
if and only if
\Eq{xyt}{
  0\leq \max\big(tf(x)+(1-t)f(y),tg(x)+(1-t)g(y)\big)
  \qquad(x,y\in X,\,t\in[0,1]).
}}

\section{Generalized convexity}

The general convexity property that we introduce below is going to play an important role in the sequel.

\Dfn{2}{Let $X$ be a nonempty set, $\circ:X\times X\to X$ be a binary operation, $p,q>0$ be constants. A function $f:X\to\R$ is called \emph{$(\circ,p,q)$-convex} if 
\Eq{*}{
  f(x\circ y)\leq pf(x)+qf(y)\qquad (x,y\in X).
}}

Trivially, if $X$ is a convex subset of a linear space, $p=q=\frac12$, and $x\circ y=\frac{x+y}{2}$, then $f$ is $(\circ,p,q)$-convex if and only if $f$ is Jensen convex.
On the other hand, if $X$ is an Abelian semigroup, $p=q=1$, and $x\circ y=x+y$, then $f$ is $(\circ,p,q)$-convex if and only if $f$ is subadditive. \\[-3mm]

The proof of the following assertion is elementary, therefore it is omitted.

\Thm{GC}{The family of $(\circ,p,q)$-convex functions is closed with respect to addition, multiplication by positive scalars and pointwise maximum.}

The main result of this paper is stated in the following theorem.

\Thm{MTGCD}{Let $X$ be a nonempty set, $\circ:X\times X\to X$ be a binary operation, and $p,q>0$ be constants. Let $f_1,\dots,f_n:X\to\R$ 
be $(\circ,p,q)$-convex functions such that
\Eq{*}{
  0\leq\max(f_1(x),\dots,f_n(x)) \qquad (x\in X).
}
Then there exist $\lambda_1,\dots,\lambda_n\geq0$ with
$\lambda_1+\cdots+\lambda_n=1$ such that
\Eq{*}{
  0\leq \lambda_1f_1(x)+\cdots+\lambda_nf_n(x) \qquad (x\in X).
}}

The follwing auxiliary result establishes the key tool for the proof of \thm{MTGCD}.

\Lem{CGC}{Let $X$ be a nonempty set, $\circ:X\times X\to X$ be a binary operation, and $p,q>0$ be constants. Let
\Eq{*}{
  S:=\Big\{\frac{a}{a+b}\,\Big|\,&\hbox{There is an operation $*:X\times X\to X$ such that} \\[-3mm]
       &\hbox{every $(\circ,p,q)$-convex function is $(*,a,b)$-convex.}\Big\}
}
Then $1-S\subseteq S$ and $S$ is dense multiplicative subsemigroup of $[0,1]$.}

\begin{proof}
If $s\in S$, then there exists an operation $*:X\times X\to X$ and $a,b>0$ such that
$s=\frac{a}{a+b}$ and $f$ is $(*,a,b)$-convex, i.e., 
\Eq{*}{
  f(x*y)\leq af(x)+bf(y)\qquad(x,y\in X).
}
Thus, interchanging the roles of $x$ and $y$, we get
\Eq{*}{
  f(y*x)\leq bf(x)+af(y)\qquad(x,y\in X),
}
which means that $f$ is $(*',b,a)$-convex, where $x*'y:=y*x$. Therefore $1-s=\frac{b}{a+b}\in S$, which shows that $1-S\subseteq S$.

Additionally, let $t\in S$ be arbitrary. Then there exists a binary operation $\cdot:X\times X\to X$ and $c,d>0$ such that $t=\frac{c}{c+d}$ and $f$ is also $(\cdot,c,d)$-convex, i.e., 
\Eq{*}{
  f(x\cdot y)\leq cf(x)+df(y)\qquad(x,y\in X).
}
Using the $(\cdot,c,d)$- and the $(*,a,b)$-convexity of $f$ (twice), for all $x,y\in X$, we obtain
\Eq{*}{
  f((x*y)\cdot (y*y))&\leq cf(x*y)+df(y*y)\\
   &\leq c(af(x)+bf(y))+d(af(y)+bf(y))\\
   &=acf(x)+(bc+ad+bd)f(y).
}
This implies that $f$ is $(\diamond,ac,bc+ad+bd)$-convex, where $x\diamond y:=(x*y)\cdot (y*y)$. Therefore, $st=\frac{ac}{ac+bc+ad+bd}\in S$, which proves that $S$ is closed with respect to multiplication.

By induction, it follows that 
\Eq{sn}{
  s^n\in S\qquad(s\in S,\,n\in\N).
}
The assumption that $f$ is $(\circ,p,q)$-convex implies that $S\cap\,]0,1[\,\neq\emptyset$. Therefore, \eq{sn} yields that $\inf S=0$. Using the inclusion $1-S\subseteq S$, we can see that $\sup S=1$. 

Finally, to prove the density of $S$ in $[0,1]$, let $0<a<b<1$ be arbitrary. By $\sup S=1$, we can choose $s\in S$ so that $\dfrac{a}{b}<s<1$. Then, for some $n\in\N$, (in particular, with $n:=\big\lfloor\frac{\log(a)}{\log(s)}\big\rfloor$), we have $s^n\in[a,b]$, which implies that $S\cap[a,b]$ is nonempty.
\end{proof}

In the next result, we verify the Maximum Theorem for two functions.

\Thm{MTTF}{Let $X$ be a nonempty set, $\circ:X\times X\to X$ be a binary operation, and $p,q>0$ be constants. If $f,g:X\to\R$ are $(\circ,p,q)$-convex functions satisfying
\Eq{MI}{
  0\leq\max(f(x),g(x)) \qquad (x\in X),
}
then there exists $\lambda\in[0,1]$ such that \eq{fgx} holds true.}

\begin{proof} First we show that $f$ and $g$ satisfy the inequality \eq{xyt}. To verify this, let $x,y\in X$ and let $s\in S$ (where the set $S$ was defined in \lem{CGC}.) Then there exist a binary operation $*:X\times X\to X$ and constans $a,b>0$ such that the $(\circ,p,q)$-convexity of $f$ and $g$ implies the $(*,a,b)$-convexity of them. Thus, by the maximum inequality \eq{MI} at $x*y$, we get
\Eq{*}{
  0\leq\max(f(x*y),g(x*y))\leq\max(af(x)+bf(y),ag(x)+bg(y)).
}
Therefore
\Eq{*}{
  0\leq\max\Big(\frac{a}{a+b}f(x)+\frac{b}{a+b}f(y),\frac{a}{a+b}g(x)+\frac{b}{a+b}g(y)\Big),
}
and hence
\Eq{*}{
  0\leq\max\big(sf(x)+(1-s)f(y),sg(x)+(1-s)g(y)\big).
}
Because $s\in S$ was arbitrary and $S$ is dense in $[0,1]$ (according to \lem{CGC}), we can conclude that \eq{xyt} is satisfied for all $t\in[0,1]$.

Having proved that \eq{xyt} is valid, in view of \cor{2F}, it follows that there exists $\lambda\in[0,1]$ such that \eq{fgx} holds.
\end{proof}

\begin{proof}[Proof of the discrete Maximum Theorem]
The statement is trivial for $n=1$ and it has been proved for $n=2$. Assume its validity for some $n\geq2$. Let $f_0,f_1,\dots,f_n$ be $(\circ,p,q)$-convex functions such that
\Eq{*}{
   0\leq\max(f_0(x),f_1(x),\dots,f_n(x))\qquad(x\in X).
}
Let $g(x):=\max(f_1(x),\dots,f_n(x))$. Then, by \thm{GC}, we have that $g$ is $(\circ,p,q)$-convex and
\Eq{*}{
   0\leq\max(f_0(x),g(x))\qquad (x\in X).
}
Using now \thm{MTTF}, we obtain the existence of $\lambda\in[0,1]$ such that
\Eq{*}{
  0&\leq \lambda f_0(x)+(1-\lambda)g(x)\\
   &=\max\big(\lambda f_0(x)+(1-\lambda)f_1(x),\dots,\lambda f_0(x)+(1-\lambda)f_n(x)\big)\qquad(x\in X).
}
By the inductive assumption, there exists $(\lambda_1,\dots,\lambda_n)\in S_n$ such that
\Eq{*}{
  0&\leq\lambda_1\big(\lambda f_0(x)+(1-\lambda)f_1(x)\big)+\cdots
    +\lambda_n\big(\lambda f_0(x)+(1-\lambda)f_n(x)\big)\\
    &=\lambda f_0(x)+\lambda_1(1-\lambda)f_1(x)+\cdots+\lambda_n(1-\lambda)f_n(x) \qquad(x\in X),
}
which proves the statement for $(n+1)$ functions. 
\end{proof}

\section{An application}

In the subsequent result we establish an extension of the 
Karush--Kuhn--Tucker Theorem.

\Thm{*}{Let $X$ be a nonempty set, $\circ:X\times X\to X$ be a binary operation, and $p,q>0$ be constants. Let $f_0,f_1,\dots,f_n:X\to\R$ be $(\circ,p,q)$-convex functions and assume that $f_0(x_0)=0$ and $x_0\in X$ is a solution of the constrained optimization problem 
\Eq{MP}{
  \mbox{Minimize }\quad f_0(x)\quad \mbox{ subject to }\quad f_1(x),\dots,f_n(x)\leq0.
}
Then there exist $(\lambda_0,\lambda_1,\dots,\lambda_n)\in S_{n+1}$ such that 
\Eq{TR}{
  \lambda_1f_1(x_0)=\dots=\lambda_1f_1(x_0)=0
}
and
\Eq{EL}{
  0\leq \lambda_0f_0(x)+\lambda_1f_1(x)+\cdots+\lambda_nf_n(x) \qquad (x\in X).
}
Conversely, if conditions \eq{TR} and \eq{EL} hold for some 
$(\lambda_0,\lambda_1,\dots,\lambda_n)\in S_{n+1}$ with $\lambda_0>0$, then $x_0$ is a solution of the optimization problem \eq{MP}.}

\begin{proof}
If $x_0$ is a solution of the optimization problem then, for all $x\in X$, the inequalities
\Eq{*}{
  f_0(x)<f_0(x_0)=0 \qquad\mbox{and}\qquad
  f_1(x),\dots,f_n(x)\leq0
}
cannot hold simultaneously. Hence
\Eq{*}{
  0\leq \max(f_0(x),f_1(x),\dots,f_n(x))\qquad(x\in X).
}
Therefore, in view of \thm{MTGCD}, there exist $(\lambda_0,\lambda_1,\dots,\lambda_n)\in S_{n+1}$ such that \eq{EL} holds.

Being a solution to \eq{MP}, $x_0$ is admissible for the optimization problem, that is, we have that $f_1(x_0),\dots,f_n(x_0)\leq0$. Hence
\Eq{*}{
  0\leq\lambda_0f_0(x_0)+\lambda_1f_1(x_0)+\cdots+\lambda_nf_n(x_0)
  =\lambda_1f_1(x_0)+\cdots+\lambda_nf_n(x_0)\leq0.
}
The terms in the last sum are nonpositive, therefore, the only way this sum can be zero is that it is zero termwise. Hence the transversality condition \eq{TR} is also true.

To prove the reversed statement, assume that \eq{TR} and \eq{EL} hold for some 
$(\lambda_0,\lambda_1,\dots,\lambda_n)\in S_{n+1}$ with $\lambda_0>0$. Let $x\in X$ be an admissible point with respect to problem \eq{MP}, i.e., assume that $f_1(x_0),\dots,f_n(x_0)\leq0$. Then, by \eq{TR} and \eq{EL}, we get
\Eq{*}{
  \lambda_0f_0(x_0)
  &=\lambda_0f_0(x_0)+\lambda_1f_1(x_0)+\cdots+\lambda_nf_n(x_0)\\
  &=0\leq\lambda_0f_0(x)+\lambda_1f_1(x)+\cdots+\lambda_nf_n(x)
  \leq \lambda_0f_0(x),
}
which, using that $\lambda_0>0$, implies $f_0(x_0)\leq f_0(x)$, and proves the minimality of $x_0$.
\end{proof}


\providecommand{\bysame}{\leavevmode\hbox to3em{\hrulefill}\thinspace}
\providecommand{\MR}{\relax\ifhmode\unskip\space\fi MR }
\providecommand{\MRhref}[2]{%
  \href{http://www.ams.org/mathscinet-getitem?mr=#1}{#2}
}
\providecommand{\href}[2]{#2}

\end{document}